\documentclass{article}
\usepackage{amsmath,amsfonts,amssymb,amsthm,graphicx}
 \usepackage{dsfont}    
\usepackage{enumerate}
\usepackage{color}
\usepackage{natbib}
\usepackage{multirow}
\usepackage{comment}

\newcommand{\fM}{M}

\usepackage{booktabs,subcaption,dcolumn}

\usepackage{algorithm}
\usepackage[noend]{algpseudocode}  
\usepackage[colorlinks=true,citecolor=blue,pdfpagemode=UseNone,pdfstartview=FitH]{hyperref}

\newif\ifFULL
\ifFULL

\fi

\newif\ifRuodu
\ifRuodu
  \renewcommand{\ge}{\geqslant}
  \renewcommand{\le}{\leqslant}
  \renewcommand{\epsilon}{\varepsilon}
  
\fi

\renewcommand{\d}{\,\mathrm{d}}

\newcommand{\p}{\mathbb{P}}

\newcommand{\E}{\mathbb{E}}    
\newcommand{\R}{\mathbb{R}}    
\newcommand{\N}{\mathbb{N}}    

  \newcommand{\id}{\mathds{1}}

\usepackage{setspace}

\theoremstyle{plain}
\newtheorem{theorem}{Theorem} 

\newtheorem{lemma}[theorem]{Lemma}

\theoremstyle{definition}

\theoremstyle{remark}
\newtheorem{remark}[theorem]{Remark}

 \renewcommand{\cite}{\citet}  

\title{The only admissible way  of merging arbitrary e-values} 
\author{
  Ruodu Wang\thanks%
  {Department of Statistics and Actuarial Science,
  University of Waterloo,
  Waterloo, Ontario, Canada.
  E-mail: \href{mailto:wang@uwaterloo.ca}{wang@uwaterloo.ca}.}}
 \date{\today}

\begin{document} 
\maketitle

\begin{abstract}
We prove that the only admissible way of merging arbitrary e-values is to use a weighted arithmetic average.
This result completes the picture of merging methods for arbitrary e-values, and generalizes the result of Vovk and Wang (2021, Annals of Statistics, 49(3), 1736--1754) that the only admissible way of symmetrically merging e-values is to use the arithmetic average combined with a constant. 
Although the proved statement is naturally anticipated, its proof relies on a sophisticated application of optimal transport duality  and a minimax theorem. 
 \vspace{2ex}
\end{abstract}

   \section{Introduction}
E-values and e-processes in hypothesis testing have been an active research area in the recent years, led by a series of papers including \cite{S21}, \cite{VW21}, \cite{GDK24} and \cite{WRB20}.  
The popularity of e-values came from its various features such as flexibility, anytime validity, robustness, post-hoc  decision validity,  and a strong connection to betting and martingales. \cite{WR22}, \cite{G24},  and \cite{RGVS23} offered many discussions on these features.

Arguably, one of the central advantages of e-values, especially in contrast to p-values, is that they are easy to combine.
It is elementary that a weighted arithmetic average of e-values, with weights summing to $1$, is an e-value.
This convenient feature has been used in many applications. 
For instance, it is used to obtain admissible ways of merging p-values by \cite{VWW22},
to construct discovery matrices by \cite{VW23},  
 to deranomize the Model-X knockoffs of \cite{BC15}
  by \cite{RB24}, and to boost the power of the e-Benjamini-Hochberg procedure of \cite{WR22} by  \cite{LR24}.

A natural question is, in addition to the weighted averages, whether there are other useful e-merging functions, i.e., functions that produce an e-value from several input e-values. 
If we additionally assume that the input e-values are independent or sequential, then many other admissible methods exist, as studied by \cite{VW21, VW24}.
Without assuming particular dependence structures, the simple question of identifying useful e-merging functions turns out to be surprisingly nontrivial. 
As a central result on the combination of e-values,
 \citet[Theorem 3.2]{VW21} showed, through several additional technical results, that all admissible symmetric e-merging functions 
take the form of a convex combination  of 
the   arithmetic average and the constant $1$.
From there, one may naturally conjecture that, without assuming symmetry,
admissible e-merging functions
should take the form of  a  convex combination of 
a weighted average and the constant $1$.
This was not established by \cite{VW21} or its follow-up papers due to technical challenges. 

The main aim of this short paper is to prove the above conjecture, and thus to offer a complete characterization of all admissible e-merging functions. 
This justifies the use of weighted averages of e-values in many applications (where the weighted averages may have already  been used), and frees us of any doubts about whether there are better choices in particular contexts. 
Although the anticipated mathematical result is natural, the proof requires a delicate analysis via optimal transport duality, very different from the proof for the symmetric case   by  \cite{VW21}. 

 At an abstract level, the reason why optimal transport duality is essential is due to 
  the requirement of an e-merging function to 
 work for e-values with any dependence structure. Maximization over all dependence structures is a classic problem in  optimal transport theory, and, more precisely,   in our context it is a multi-marginal optimal transport problem; see \cite{RR98} and \cite{V09} for more on optimal transport theory.  On a related matter, optimal transport duality is crucial for results of \cite{VWW22}, where e-values serve as intermediate tools for characterizing admissible ways of merging p-values.

Merging functions of some subclasses of e-values are discussed in Section \ref{sec:r1-5}.
 
  \section{The main result}

 Let $K$ be a positive integer and denote by $ [K]=\{1,\dots,K\}$.
Write $\R_+= [0,\infty)$ and denote by $\Delta_n$ the standard simplex in $\R^n$  for  $n\in \N$, that is,  $$\Delta_n=\left\{(x_1,\dots,x_n)\in [0,1]^n: \sum_{k=1}^n x_k=1 \right\}.$$  
A \emph{hypothesis} is a collection of probability measures on a measurable space $(\Omega,\mathcal F)$  called the sample space.

 An \emph{e-variable} $E$  for a hypothesis $\mathcal H$ is a $[0,\infty]$-valued random variable satisfying $\E^\p [E]\le1$ for all $\p \in \mathcal H$.   
 We have been using the loose term ``e-values" for e-variables in the Introduction.
  An \emph{e-merging function}    is an increasing  Borel function
$F:\R_+^K\to\R_+$
such that for any hypothesis, 
$
  F(E_1,\dots,E_K) $ is an e-variable for any e-variables $E_1,\dots,E_K$. We have safely excluded infinite points in the domain of $F$; see \citet[Appendix C]{VW21} for a justification. 
  
  An e-merging function $F$ is \emph{admissible} if 
for any e-merging function $G$, $G\ge F$ implies $G=F$.
That is,  $F$ cannot be strictly improved.
All equalities and inequalities are point-wise.

 As explained in \citet[Appendix D]{VW21}, to study  merging methods, it is necessary and sufficient to consider the simple hypothesis $\{\p\}$ with an atomless probability measure $\p$, 
 and all results carry through to any settings of possibly composite hypotheses. Therefore, in the sequel, all e-variables are defined for the fixed $\mathcal H=\{\p\}$, and  $\E$ represents the expectation with respect to $\p$. We omit mentioning $\p$ explicitly in most places.

For any vector $\boldsymbol \lambda \in \Delta_{K+1}$, 
define the mapping $\fM_{\boldsymbol \lambda}: \R_+^{K}\to \R_+$ 
by $\mathbf e\mapsto   \boldsymbol \lambda \cdot (\mathbf e,1)$.
The mapping $\fM_{\boldsymbol \lambda}$ outputs the $ {\boldsymbol \lambda}$-weighted average of its input arguments and the constant $1$. 
It is an e-merging function of dimension $K$. 
The main result in the paper is stated next.

 \begin{theorem}
 \label{th:admissible}
For a function $ F:\R_+^K\to\R_+$,
 \begin{enumerate}[(i)]
 \item if $F$ is an    e-merging function, then
  $F\le \fM_{ \boldsymbol \lambda }$ for some  $\boldsymbol \lambda \in \Delta_{K+1}$;  
\item $F$ is an  admissible e-merging function  if and only if     $F = \fM_{ \boldsymbol \lambda }$ for some $\boldsymbol \lambda \in \Delta_{K+1}$. 
\end{enumerate}
\end{theorem}

It suffices to show (i), as (ii) follows directly from (i) and the simple fact that $ \fM_{ \boldsymbol \lambda }$ is an e-merging function not dominated by any other weighted average. 
In Section \ref{sec:pf} we will prove Theorem \ref{th:admissible}.
Before that, we address the simplest case $K=1$, which will be used in the proof of Theorem \ref{th:admissible}.
In what follows, a random variable is binary if it   takes at most two possible values. 
\begin{lemma}\label{lem:mean2}
Let $\theta\in[ 1,\infty]$ and $r\in \R_+$. 
If a function $g:[0,\theta] \to \R_+$ satisfies 
$\E[g(E)]\le r  $ for all binary e-variables $E$ taking     values in $[0,\theta]$, then there exists $h \in [0,1]$ such that  $g(x)\le r (1-h + h x )$ for $x\in [0,\theta]\cap \R$.
\end{lemma}


As a particular case of Lemma \ref{lem:mean2} with $r=1$, if $g$ is an e-merging function of dimension $1$, then  there exists $h \in [0,1]$ such that 
$g(x) \le (1-h) + h x  =\fM_{(h,1-h)}(x)$ for $x\in \R_+$.
Hence, Theorem \ref{th:admissible} holds for $K=1$.
This result, although very simple, is   not found in the literature.

Lemma \ref{lem:mean2} is a refinement of Lemma EC.2 of \cite{WWZ22}; the latter result considered all e-variables instead of those binary and bounded in $[0,\theta]$. 
Because of its importance for our main result, 
a self-contained proof,  similar to the one in \cite{WWZ22}, is presented in Section \ref{sec:lem2}.
 
\begin{remark}
For $\alpha \in [0,1]$, a level-$\alpha$ test for a hypothesis $\mathcal H$ is a $[0,1]$-valued random variable $\tau$ satisfying $\E^\p [\tau ]\le \alpha$ for all $\p \in \mathcal H$,
with $\tau=1$ indicating a rejection, $\tau=0$ indicating no rejection, and $\tau\in (0,1)$ indicating a randomized decision. 
For $\alpha \ne 0$, $\tau/\alpha$ is an e-variable. For the connection between tests and e-values, 
see \cite{K24} and \cite{RW24}.
Applying the proof of Theorem \ref{th:admissible} to  e-variables that are bounded in $[0,1/\alpha]$, 
we get that using a weighted average is the only admissible way to merge arbitrary level-$\alpha$ tests.
Merging tests at different levels is similar, up to an adjustment of the weights.
\end{remark}
  \section{Proof of Theorem \ref{th:admissible}}

\label{sec:pf}

Let us first prove a small lemma. In what follows, $a \vee b$ is the maximum of $a,b$, $a \wedge b$ is the minimum of $a,b$,
and $\mathbf 0$ represents a zero vector of the appropriate dimension.
\begin{lemma}\label{lem:dom}
Any  e-merging function $ F $ satisfies  $F( \mathbf e)\le  1 \vee \max( \mathbf e)$ for all $ \mathbf e\in \R_+^K$.
\end{lemma}
\begin{proof}
Suppose that there exists $ \mathbf e \in \R_+^K$ such that 
$F( \mathbf e ) > 1 \vee \max( \mathbf e)$.
Let  $\bar e=\max \{ \mathbf e\}$.
If $ \bar e \le 1$, then $ \mathbf e$ is a vector of constant e-variables,
but $F( \mathbf e)>1$ is not an e-variable, a contradiction. 
Next, consider $ \bar e > 1$.
Take e-variables $E_1,\dots,E_K$ with 
$\p((E_1,\dots,E_K)= \mathbf e)=1/\bar e$ 
and $\p((E_1,\dots,E_K)=\mathbf 0)=1-1/\bar e$.
It is straightforward to see that $E_1,\dots,E_K$ are e-variables.
Moreover, $\E[F(E_1,\dots,E_K)]> \bar e /\bar e=1$, a contradiction to the assumption that $F$ is an e-merging function. 
\end{proof}

Fix an e-merging function $F$. 
Lemma \ref{lem:dom} implies in particular  
\begin{align}\label{eq:dominated}
F(e_1,\dots,e_K)\le 1+ \sum_{k=1}^K e_k
\mbox{ ~~~ for all $(e_1,\dots,e_K)\in \R_+^K$.}
\end{align}
This allows us to apply optimal transport duality.
Let us state the duality first. 
Let $\Pi(\mu_1,\dots,\mu_K)$ be the set of all Borel measures on $\R^K$ with marginal distributions $\mu_1,\dots,\mu_K$ on $\R$.  
Let $\mathcal B$ be the set of Borel-measurable functions on $\R_+$, and  
write $\bigoplus_{k=1}^K \phi_k :(x_1,\dots,x_K)\mapsto  \sum_{k=1}^K \phi_k (x_k)$ for functions $\phi_1,\dots,\phi_K\in \mathcal B $.
Define  
$$D_F=\left\{(\phi_1,\dots,\phi_K) \in \mathcal B^K: \bigoplus_{k=1}^K   \phi_k   \ge F\right\},$$
that is, the set of tuples of univariate functions  whose sum  dominates $F$.  
The set $D_F$ is important in optimal transport theory because it compares $F$ with a \emph{separable} function $G$ of the from $G=\bigoplus_{k=1}^K   \phi_k $, with the key property that $\E[G(E_1,\dots,E_K )]$  depends only on the marginals of $(E_1,\dots,E_K)$ but not the dependence structure.

 Let $\mathcal M_{\mathcal E}$ be the set of all distributions on $\R_+$ with mean no larger than $1$, i.e., the set of all distributions of e-variables,
and   $\mathcal M_{\mathcal E}^\theta $ be the subset of $\mathcal M_{\mathcal E} $ containing all distributions on $[0,\theta ]$ for $\theta \ge 1$.
In what follows, we always write $\boldsymbol \phi=(\phi_1,\dots,\phi_K)$
and 
$\boldsymbol \mu=(\mu_1,\dots,\mu_K)$.

Condition \eqref{eq:dominated} implies that the set $D_F$ is not empty, and moreover some elements $\boldsymbol \phi$ in $D_F$ satisfies  $\sum_{k=1}^K\int \phi_k \d \mu_k<\infty$ when $\boldsymbol \mu \in (\mathcal M_{\mathcal E})^K$. 
With this, optimal transport duality holds. 

\begin{lemma}[Optimal transport duality] 
\label{lem:OTD}
For $\mu_1,\dots,\mu_K\in \mathcal  M_{\mathcal E}$ and an e-merging function $F$, we have
\begin{align}
\label{eq:dual}
 1\ge \sup_{\pi\in \Pi(\boldsymbol \mu) } \int F \d \pi = \inf_{\boldsymbol \phi \in D_F}  \sum_{k=1}^K\int \phi_k \d \mu_k,
\end{align}
where in the infimum we only consider those with $  \sum_{k=1}^K\int \phi_k \d \mu_k$ well-defined.
\end{lemma}
The inequality in \eqref{eq:dual} follows because $F$ is an e-merging function, and the equality is a classic form of optimal transport duality in the form of Remark 2.1.2 of \cite{RR98}.

To interpret \eqref{eq:dual}, 
let $E_1,\dots,E_K$ be e-variables  with distributions $\mu_1,\dots,\mu_K$. 
Since $F \le  \bigoplus_{k=1}^K \phi_k $ for $\boldsymbol \phi \in D_F$, the equality in \eqref{eq:dual} means that 
 the supremum of $\E[F(E_1,\dots,E_K )]  $ over dependence structures of $(E_1,\dots,E_K)$ is equal to the infimum of $\E[ \sum_{k=1}^K \phi_k(E_k)]$ over $\boldsymbol \phi $.

We explain our main idea in the proof first.
The inequality \eqref{eq:dual}  
implies 
$$ 
\sup_{\boldsymbol \mu \in (\mathcal M_{\mathcal E})^K} \inf_{\boldsymbol \phi \in D_F}  \sum_{k=1}^K\int \phi_k \d \mu_k \le 1. 
$$ 
If we could interchange the order of sup and inf in the above equation,  
 then for any $\epsilon>0$,  
$  \sum_{k=1}^K \int \phi_k\d \mu_k \le 1+\epsilon$
for some $\boldsymbol\phi$ and all $\boldsymbol \mu \in (\mathcal M_{\mathcal E})^K$. 
Using this and Lemma \ref{lem:mean2} 
 yields  linear upper bounds  on   $\boldsymbol\phi$, with which we can further use $F\le \bigoplus_{k=1}^K   \phi_k $ to show   $F\le M_{\boldsymbol \lambda}$ for some $\boldsymbol \lambda$.
 The following lemma gives  the desired maximin interchangeability under the assumption of bounded compact.
  Let $ \mathcal M_0 $ be the set of all distributions on $\R_+$.

\begin{lemma}
\label{lem:sion2}
Suppose   $\mathcal M  \subseteq ( \mathcal M_0)^K$ is compact with respect to weak convergence  and $F:\R_+^K\to \R_+$ is bounded and upper semicontinuous. Then 
\begin{equation}
\label{eq:sion2}
  \sup_{ \boldsymbol\mu \in   \mathcal M } \inf_{\boldsymbol \phi \in D_F}  \sum_{k=1}^K\int \phi_k \d \mu_k  
  = \inf_{\boldsymbol \phi \in D^+_F}    \sup_{ \boldsymbol\mu \in   \mathcal M }  \sum_{k=1}^K\int \phi_k \d \mu_k, 
\end{equation} 
where $D_F^+$ is the subset of $D_F$ containing all $\boldsymbol \phi$ with nonnegative components. 
\end{lemma}

We will prove Lemma \ref{lem:sion2} later. We now use this lemma to complete the proof of Theorem \ref{th:admissible}. 
Fix $\theta \in[ 1,\infty)$. We need to verify a few things for an e-merging function $F$, which justifies that we can apply Lemma \ref{lem:sion2}.
\begin{enumerate}[(i)]
\item 
Since $[0,\theta ]$ is compact, the set $(\mathcal M ^\theta_{\mathcal E})^K$ equipped with the weak topology is tight.
   By Prokhorov's theorem, it is   sequentially compact, and hence compact. 
   \item 
Since $F$ is bounded on $[0,\theta]^K$ and its value outside $[0,\theta]^K$ does not matter in \eqref{eq:sion2}, we can treat $F$ as bounded when applying Lemma \ref{lem:sion2}. 
\item There exists an upper semicontinuous e-merging function $F^*$ with $F^*\ge F$   by Proposition E.1 of \cite{VW21}. Therefore, to prove Theorem \ref{th:admissible}, it suffices to consider upper semicontinuous e-merging functions. 
\end{enumerate} 

 By   using
Lemmas \ref{lem:OTD} and \ref{lem:sion2} with $\mathcal M= (\mathcal M ^\theta_{\mathcal E})^K$  and  an upper semicontinuous e-merging function $F$, we have   
\begin{align}
1& \ge \inf_{\boldsymbol \phi \in D^+_F}  \sup_{\boldsymbol \mu \in (\mathcal M^\theta _{\mathcal E})^K}  \sum_{k=1}^K\int \phi_k \d \mu_k  = \inf_{ \boldsymbol \phi  \in D^+_F} \sum_{k=1}^K  \sup_{ \mu_k  \in \mathcal M^\theta_{\mathcal E}}   \int \phi_k \d \mu_k . \label{eq:value}
 \end{align}
 Denote by $T_{\phi}= \sup_{ \mu \in \mathcal M_{\mathcal E}^\theta}   \int \phi  \d \mu  $ for $\phi\in \mathcal B$. 
 For any $\epsilon>0$, by   \eqref{eq:value},  we can find  $(\phi_1,\dots,\phi_K)\in D^+ _F$
 such that $\sum_{k=1}^K T_{\phi_k}\le 1+\epsilon$.
For  each $k\in [K]$, by applying Lemma \ref{lem:mean2} to  $\phi_k$, there exists a constant $h_{\phi_k}\in [0,1]$ such that 
 $$ \phi_k(x) \le T_{\phi_k}\left(1-h_{\phi_k} + h_{\phi_k} x\right) \mbox{~for all $x\in [0,\theta]$.}$$ 
Since $(\phi_1,\dots,\phi_K)\in D^+_F$, we have
$$
F(x_1,\dots,x_K) \le\sum_{k=1}^K  T_{\phi_k}\left(1-h_{\phi_k} + h_{\phi_k} x_k\right)   \mbox{~for all  $x_1,\dots,x_K \in [0,\theta]$}.
$$
This means that there exists $\boldsymbol \lambda_{\epsilon,\theta} \in \Delta_{K+1}$
such that   $F \le (1+\epsilon) \fM_{ \boldsymbol \lambda_{\epsilon,\theta} }$ on $ [0,\theta]^K$. 
Since $\Delta_{K+1}$ is compact, 
we can find a convergent subsequence of $\boldsymbol \lambda_{\epsilon,\theta} $  by sending $\epsilon \downarrow 0$ and $\theta\to \infty$, with its limit denoted by $ \boldsymbol \lambda_0 \in \Delta_{K+1}$.
Continuity of  $ \boldsymbol \lambda 
\mapsto \fM_{ \boldsymbol \lambda  }$   
yields $ F \le  \fM_{ \boldsymbol \lambda_0  } $ on $ \R_+^K$, 
thus the desired statement. \qed

\begin{proof}[of Lemma \ref{lem:sion2}] 
We 
  assume $0\le F\le 1$ without loss of generality.   
We first argue   that we can require $(\phi_1,\dots,\phi_K)\in D_F$   to
be upper semicontinuous and
 take values in $[0,1]$.  
 That is, 
 we will show, for $(\mu_1,\dots,\mu_K )\in \mathcal M $,
 \begin{align}
\label{eq:dual3}
  \inf_{\boldsymbol \phi \in D_F}  \sum_{k=1}^K\int \phi_k \d \mu_k= \inf_{\boldsymbol \phi \in \widetilde D_F}  \sum_{k=1}^K\int \phi_k \d \mu_k,
\end{align}  
where
$$ \widetilde D_F=\left\{ \boldsymbol \phi \in  D_F:  \mbox{ $\phi_k$ is upper semicontinuous  and  $0\le  \phi_k\le  1$  for $k\in [K]$}\right\}.$$

\begin{enumerate}[(a)]
\item For  $(\phi_1,\dots,\phi_K)\in D_F$, 
 denote by $c_k= \inf_{x\in \R_+} \phi_k(x)$ for $k\in [K]$. 
 Since $\sum_{k=1}^K \phi_k(x_k) \ge F(x_1,\dots,x_K) \ge 0$ for all $(x_1,\dots,x_K)\in \R_+^K$,
 the function $\tilde\phi_k: x\mapsto \phi_k(x)-c_k +\sum_{m=1}^K c_m/K$ is nonnegative. 
 Moreover, for all $(\mu_1,\dots,\mu_K)\in \mathcal  M $,
 $$
 \sum_{k=1}^K\int \tilde \phi_k \d \mu_k =    \sum_{k=1}^K\int \phi_k \d \mu_k   .
 $$
Hence, we can safely restrict the functions   $\phi_1,\dots,\phi_K$ to be nonnegative. 
\item For  $(\mu_1,\dots,\mu_K)\in \mathcal  M $,
 since $F \le 1$, together with (a) we can also safely truncate $\phi_1,\dots,\phi_K$ at $ 1$. 
\item Since $F$ is upper semicontinuous and bounded, by Proposition 1.31 of \cite{K84} the functions $(\phi_1,\dots,\phi_K)$ can be chosen as upper semicontinuous.
\end{enumerate}
The above three arguments together show the assertion \eqref{eq:dual3}.

Next, we verify 
\begin{align}
 \sup_{\boldsymbol \mu \in  \mathcal M } \inf_{ \boldsymbol \phi\in  \widetilde D _F}  \sum_{k=1}^K  
 \int \phi_k \d \mu_k  
  = \inf_{\boldsymbol \phi \in \widetilde D _F}  \sup_{\boldsymbol \mu \in  \mathcal M }   \sum_{k=1}^K\int \phi_k \d \mu_k.
  \label{eq:inter2}
   \end{align} 
To apply the minimax theorem of \cite{S58}, we need to check a few things. 
Define the mapping $$J:( \boldsymbol \mu,\boldsymbol \phi)\mapsto \sum_{k=1}^K  
 \int \phi_k \d \mu_k .$$
 \begin{enumerate}[(i)]
 \item  
 The mapping $J$ is bilinear, and therefore both convex and concave. 
 \item The set $\mathcal M$ equipped with the weak topology is compact.
   \item 
Since each component of $ \boldsymbol \phi $ is bounded in $[0, 1]$,
$\boldsymbol \phi \mapsto J( \boldsymbol \mu,\boldsymbol \phi)$ is continuous with respect to 
 point-wise convergence in $\boldsymbol \phi$ for each $\boldsymbol \mu$. This continuity is not needed if we use the minimax theorem of \citet[Theorem 2]{F53}.
    \item 
Since each component of $ \boldsymbol \phi $ is upper semicontinuous,
$\boldsymbol \mu \mapsto J( \boldsymbol \mu,\boldsymbol \phi)$ is upper semicontinuous with respect to 
weak convergence in $\boldsymbol \mu$ for each $\boldsymbol \phi$. 
   \end{enumerate}
Therefore, Sion's minimax theorem implies that \eqref{eq:inter2} holds.  
Putting  \eqref{eq:dual3} and \eqref{eq:inter2} together we get
$$\sup_{ \boldsymbol\mu \in   \mathcal M } \inf_{\boldsymbol \phi \in D_F } J( \boldsymbol \mu,\boldsymbol \phi)   = \inf_{\boldsymbol \phi \in \widetilde D_F}  \sup_{\boldsymbol \mu \in  \mathcal M  }J( \boldsymbol \mu,\boldsymbol \phi).
$$ 
 Finally, by using $ \widetilde D_F  \subseteq D_F^+ \subseteq D_F$ we obtain  \eqref{eq:sion2}.
\end{proof}

\section{Proof of Lemma \ref{lem:mean2}}
\label{sec:lem2}
\begin{proof}
 
The case $r=0$ is trivial as $g$ is the constant function $0$. 
Otherwise, $r>0$, and without loss of generality we can assume $r=1$.

If $\theta=1$, then   $g(x)\le 1$ for all $x$, and taking $h=0$ gives the desired inequality. 
In what follows, we assume $\theta \in (1,\infty]$,
and all points $x,y$ that appear below are in $[0,\theta]\cap \R$.

 First, it is easy to note that $g(y)\le y$ for $y>1$; indeed, if $g(y)>y$ then taking a random variable $X$ with $\p(X=y)=1/y$ and $0$ otherwise  gives $\E[g(X)]> 1$ and breaks the assumption.
 Moreover, $g(y)\le 1$ for $y\le 1$ is also clear, which in particular implies $g(1)\le 1$.

Suppose for the purpose of contradiction that the statement in the lemma does not hold. 
This means that for each $\lambda\in [0,1]$, either (a) $g(x) >(1-\lambda) + \lambda x$ for some $x<1$ or (b) $g(y) >(1-\lambda) + \lambda y$ for some $y>1$ (or both).
Since $g(y)\le y$ for $y>1$ and $g(x)\le 1$ for $x<1$, we know that
$\lambda =1$ implies (a)
and  $\lambda =0$ 
implies (b). 

We claim that there exists $\lambda_0\in (0,1)$ for which both (a) and (b) happen. 
To show this claim, let 
\begin{align*} 
\Lambda_0&= \{ \lambda \in [0,1]: g(y) >(1-\lambda) + \lambda y \mbox{ for some $y>1$}\};\\\Lambda_1&= \{ \lambda \in [0,1]: g(x) >(1-\lambda) + \lambda x \mbox{ for some $x<1$}\}.\end{align*}
Clearly, the above arguments show $\Lambda_0\cup \Lambda_1=[0,1]$, $0\in \Lambda_0$, and $1\in \Lambda_1$. Moreover, since the function $\lambda 
\mapsto  (1-\lambda) + \lambda x$ is monotone  for either $x<1$ or $x>1$, we know that both $\Lambda_0$ and $\Lambda_1$ are intervals. 
Let  
$ \lambda_*=\sup \Lambda_0$ and  
$ \lambda^*=\inf \Lambda_1.$
We will argue  $\lambda_*\not \in \Lambda_0$ 
and $ \lambda^* \not \in   \Lambda_1.$ 
If $\lambda_* \in \Lambda_0$, then there exists $y>1$ such that   $g(y) >(1-\lambda_*) + \lambda_* y $. By continuity, there exists $\hat \lambda_*>\lambda_*$ such that  $g(y) >(1-\hat \lambda_*) + \hat\lambda_* y $, showing that $\hat \lambda_*\in \Lambda_0$, contradicting the definition of $\lambda_*$. Therefore, $\lambda_* \not \in \Lambda_0$.
Similarly,  $\lambda^* \not \in \Lambda_1$, following the same argument.
If $\lambda_*  = \lambda^*$, then this point is not contained in $\Lambda_0\cup \Lambda_1$, a contradiction to $\Lambda_0\cup \Lambda_1=[0,1]$. Hence, it must be  $\lambda_* > \lambda^*$, which implies that $\Lambda_0\cap \Lambda_1 $ is not empty.

Let $x_0<1$ and $y_0>1$ be such that 
  $$ g(x_0)> 1-\lambda_0 +  \lambda_0 x_0  \mbox{~~~and~~~} g(y_0)> 1-\lambda_0 +  \lambda_0 y_0 .$$    
Let $X$ be distributed as 
$\p(X=y_0) = (1-x_0)/(y_0-x_0)$
and 
 $\p(X=x_0) = ( y_0-1)/(y_0-x_0)$, which clearly satisfies 
 $\E[X]=1$ and is binary.   
Moreover, 
\begin{align*}
\E[g(X)] &=
\frac{ 1-x_0   }{y_0-x_0 }  g(y_0 ) + 
\frac{ y_0-1   }{y_0-x_0 } g(x_0)
\\&> 
\frac{ 1-x_0   }{y_0-x_0 } \left(1-\lambda_0 +  \lambda_0 y_0 \right) + 
\frac{ y_0-1   }{y_0-x_0 }\left(1-\lambda_0 +  \lambda_0 x_0\right)=1.
\end{align*}
This yields a contradiction.  
\end{proof}

\section{Merging functions of subclasses of e-values}
\label{sec:r1-5}

 Theorem \ref{th:admissible} focuses on merging arbitrary e-variables.
In an application, if some information on the e-variables 
is available, one may consider 
merging functions that are valid   for a subclass of e-variables, but not necessarily valid for arbitrary e-variables.
For many subclasses, every weighted average $\fM_{ \boldsymbol \lambda }$ for  $\boldsymbol \lambda \in \Delta_{K+1}$ remains admissible, but they may not be the only admissible choices.

Formally, let $\mathcal E_K$ be a set of $K$-dimensional vectors of e-variables.
A \emph{merging function 
of the subclass $\mathcal E_K$}  
 is an increasing  Borel function
$F:\R_+^K\to\R_+$
such that 
$
  F(E_1,\dots,E_K) $ is an e-variable for any $(E_1,\dots,E_K)\in \mathcal E_K$.
A merging function $F$ 
of $\mathcal E_K$    is \emph{admissible} if 
for any  merging function $G$ of $\mathcal E_K$, $G\ge F$ implies $G=F$.
   The cases of independent and sequential e-variables are studied by \cite{VW21,VW24}.
     We  consider four other subclasses  below.  
     
We first make a general observation. 
   Let $\mathbf E$ be a vector of e-variables with  mean $1$.
 Clearly, $\E[M_{\boldsymbol \lambda} (\mathbf E)]=1$ for $\boldsymbol \lambda\in \Delta_{K+1}$.
 Since $M_{\boldsymbol \lambda}$ is continuous, any $G\ge M_{\boldsymbol \lambda} $ with $G\ne  M_{\boldsymbol \lambda}$    satisfies 
 $\E[G(\mathbf E)]>1$ if $\mathbf E$ has full support in  $\R_+^K$, and $G$ cannot be a valid merging function of $\{\mathbf E\}$. Therefore, if the subclass $\mathcal E_K$ 
contains an element of full support in $\R_+^K$, then  $M_{\boldsymbol \lambda}$ is admissible.
 
(a) \textbf{Fixed Marginals.}
Let $\mathcal E_K$ be the set of vectors of  e-variables  with given marginals $\mu_1,\dots,\mu_K$. 
If all marginals have full support in $\R_+ $, then 
$M_{\boldsymbol \lambda}$ is  admissible.
To see other admissible merging functions of $\mathcal E_K$, 
suppose that the marginals $\mu_1,\dots,\mu_K$ are continuous with full support on $\R_+ $ and let $g_1,\dots,g_K$ be the survival functions of $\mu_1,\dots,\mu_K$. 
Let 
$f_1,\dots,f_K$ be admissible calibrators; that is, decreasing upper semicontinuous  functions $f_k:[0,1] \to [0,\infty]$ satisfying $f_k(0)=\infty$ and $\int_0^1 f_k(x) \d x=1$. 
For each $k\in [K]$ and $E_k$   distributed as $\mu_k$,  $g_k(E_k)$ is uniformly distributed on $[0,1]$, and  $f_k(g_k(E_k))$ is an e-variable. Therefore,   
 $$
  (e_1,\dots,e_K) \mapsto M_{\boldsymbol \lambda}(f_1 \circ g_1(e_1) ,\dots, f_K\circ g_K(e_K) ) 
 $$
 is an admissible merging function of $\mathcal E_K$, following the general observation above.

(b) \textbf{Bounded second moments.} Let $\mathcal E_K$ be the set of vectors $(E_1,\dots,E_K)$ of  e-variables  with second moments $\E[X_iX_j]$ bounded above by $\sigma_{ij}\in [0,\infty)$ for $i,j\in [K]$. In addition to the weighted averages, the functions $(e_1,\dots,e_K)\mapsto  e_ie_j/\sigma_{ij}$ and their mixtures are merging functions of $\mathcal E_K$.  
As a special case, if elements of $\mathcal E_K$ have nonpositive bivariate correlation coefficients (a form of negative dependence), then the functions  $(e_1,\dots,e_K)\mapsto  e_ie_j $ with $i\ne j$ are merging functions  of $\mathcal E_K$.  

(c) \textbf{Identical e-variables.}
Let $\mathcal E_K$ be the set of vectors of identical e-variables. 
Clearly, the function  $V_K:\mathbf e \mapsto \max (\mathbf e)$ is a merging function of the subclass $\mathcal E_K$.  Moreover, 
all admissible merging functions in this setting have the form $\lambda +(1-\lambda)V_K$ for $\lambda \in[0,1]$.
To see this, first note that for any merging function $F$ of   $\mathcal E_K$, the function $g:e\mapsto F(e,\dots,e)$ satisfies the condition in Lemma \ref{lem:mean2},
and hence for some $\lambda \in[0,1]$ $g(e)\le \lambda + (1-\lambda) e$ for all $e\in \R_+$. 
It follows that $F(\mathbf e) \le g(\max (\mathbf e)) \le \lambda + (1-\lambda) \max (\mathbf e)$ for all $\mathbf e\in \R_+^K$.
Note that when applied to a vector $\mathbf E$ of identical e-variables,
 $V_K(\mathbf E)=M_{\boldsymbol \lambda}(\mathbf E)$ for any $\boldsymbol \lambda$ with the final component $0$, so using $ \lambda + (1-\lambda)V_K$ is the same as using weighted averages. 

(d) \textbf{Exchangeable e-variables.} 
Let $\mathcal E_K$ be the set of vectors of exchangeable e-variables. 
For $(E_1,\dots,E_K)\in \mathcal E_K$,  let $\bar E_k=(E_1+\dots+E_k)/k$ for $k\in [K]$.
By the randomized Markov inequality in   \citet[Theorem 4.5]{RW24},
$$
\E \left [\id_{\{\bigvee_{k\in [K]} \bar E_k\ge \beta \}} \right] \le 1/\beta
\mbox{~~~for all $\beta>1$.}
$$ 
Hence,  $F_\beta:(e_1,\dots,e_K)\mapsto  \beta\id_{\{\bigvee_{k\in [K]} \bar e_k\ge \beta \}}$ for $\beta>1$ is an e-merging function of $\mathcal E_K$.
It is straightforward to check that neither $F_\beta \ge M_{\boldsymbol \lambda}$ nor $F_\beta\le M_{\boldsymbol \lambda}$ for any $\beta>0$ and $\boldsymbol \lambda\in \Delta_{K+1}$. 
It remains unclear whether $F_\beta $ is admissible.

\subsection*{Acknowledgements}

The author would like to thank the editors,  two anonymous referees for useful comments on the paper, and
 Jose Blanchet, Qinyu Wu, Vladimir Vovk, and Zhenyuan Zhang for helpful discussions on the main result.
The work is financially supported by the Natural Sciences and Engineering Research Council of Canada (RGPIN-2024-03728, CRC-2022-00141).


\end{document}